\documentclass{ifacconf}

\usepackage{graphicx}      
\usepackage{natbib}        

\usepackage{amsmath}
\usepackage{mathtools}
\usepackage{amssymb}
\newcommand{\ie}{i.\,e.~}
\newcommand{\eg}{e.\,g.~}
\newcommand{\Span}[1]{\mathrm{span}\{#1\}}
\newcommand{\Rank}[1]{\mathrm{rank}\left(#1\right)}
\newcommand{\Dim}[1]{\mathrm{dim}(#1)}

\newcommand{\D}{\mathrm{d}}
\newcommand{\Lie}{\mathrm{L}}

\newcommand{\Cos}{\mathrm{c}}
\newcommand{\Sin}{\mathrm{s}}
\newcommand{\Tan}{\mathrm{t}}

\begin{document}
	\begin{frontmatter}
	
		\title{On the Linearization of Flat Two-Input Systems by Prolongations and Applications to Control Design} 
		
		\thanks[footnoteinfo]{The first author and the second author have been supported by the Austrian Science Fund (FWF) under grant number P 32151 and P 29964.}
		
		\author[First]{Conrad Gst{\"o}ttner} 
		\author[First]{Bernd Kolar} 
		\author[First]{Markus Sch{\"o}berl}
		
		\address[First]{Institute of Automatic Control and Control Systems Technology, Johannes Kepler University, Linz, Austria.\\(e-mail: \{conrad.gstoettner,bernd.kolar,markus.schoeberl\}@jku.at)}
		
		\begin{abstract}
			In this paper we consider $(x,u)$-flat nonlinear control systems with two inputs, and show that every such system can be rendered static feedback linearizable by prolongations of a suitably chosen control. This result is not only of theoretical interest, but has also important implications on the design of flatness based tracking controls. We show that a tracking control based on quasi-static state feedback can always be designed in such a way that only measurements of a (classical) state of the system, and not measurements of a generalized Brunovsky state, as reported in the literature, are required.
		\end{abstract}
		
	
	\end{frontmatter}
	
	\section{Introduction}
		The concept of flatness was introduced in control theory by Fliess, L\'evine, Martin and Rouchon, see \eg \cite{FliessLevineMartinRouchon:1992,FliessLevineMartinRouchon:1995}, and has attracted a lot of interest in the control systems theory community. The flatness property allows an elegant systematic solution of feed-forward and feedback problems, see \eg \cite{FliessLevineMartinRouchon:1995}. Roughly speaking, a nonlinear control system
		\begin{align}\label{eq:1}
			\begin{aligned}
				\dot{x}&=f(x,u)
			\end{aligned}
		\end{align}
		with $\Dim{x}=n$ states and $\Dim{u}=m$ inputs is flat, if there exist $m$ differentially independent functions $y^j=\varphi^j(x,u,u_1,\ldots,u_q)$, $u_k$ denoting the $k$-th time derivative of $u$, such that $x$ and $u$ can be parameterized by $y$ and its time derivatives. For this parameterization we write
		\begin{align*}
			\begin{aligned}
				x&=F_x(y,y_1,\ldots,y_{r-1})\\
				u&=F_u(y,y_1,\ldots,y_r)
			\end{aligned}
		\end{align*}
		and refer to it as parameterizing map with respect to the flat output $y$. For a given flat output, $F_x$ and $F_u$ are unique. If the parameterizing map is invertible, \ie $y$ and all the time derivatives of $y$ present in the map can be expressed solely as functions of $x$ and $u$, the system is static feedback linearizable. In this case we call $y$ a linearizing output of the static feedback linearizable system. In contrast to the static feedback linearization problem, which is completely solved, see \cite{JakubczykRespondek:1980}, \cite{NijmeijervanderSchaft:1990}, there are many open problems concerning flatness. Recent research in the field of flatness can be found in \eg \cite{SchlacherSchoberl:2013}, \cite{SchoberlSchlacher:2014}, \cite{NicolauRespondek:2017}.
		
		In this paper we confine ourselves to systems of the form \eqref{eq:1} with two inputs which are $(x,u)$-flat, \ie systems which possess a flat output of the form $y=\varphi(x,u)$, which may depend on $u$ but not on time derivatives of $u$. We assume that the systems we deal with have no redundant inputs, \ie $\Rank{\partial_{u}f}=m$. Furthermore, we assume that all functions we deal with are smooth.
		
		It is well known that every flat system can be rendered static feedback linearizable by an endogenous dynamic feedback. If a flat output is known, such an endogenous dynamic feedback can be constructed in a systematic way, see \eg \cite{FliessLevineMartinRouchon:1999}. In the present paper, we deal with the linearization by a special sub-class of endogenous dynamic feedbacks, namely by (repeated) prolongations of a suitable control. Systems that are linearizable by one-fold prolongation of a suitable control are considered \eg in \cite{NicolauRespondek:2017}, where a complete solution of the flatness problem for this class of systems is provided. The linearization by other restricted classes of dynamic feedbacks is considered \eg in \cite{CharletLevine:1989} and \cite{CharletLevineMarino:1991}. Furthermore, results on systems linearizable by prolongations of the original inputs can be found in \cite{SluisTilbury:1996}, \cite{FossasFranch:2000} and \cite{FranchFossas:2005}. In the present contribution, we prove that every $(x,u)$-flat system with two inputs is linearizable by prolongations of a single (new) input after a suitable static feedback has been applied.
		
		This theoretical result concerning the linearization by prolongations is very useful for the design of quasi-static flatness based tracking controls. Tracking control based on exact linearization by a quasi-static state feedback can be found in \cite{DelaleauRudolph:1998}, and has the advantage that it results in a static control law. However, it requires measurements (or estimates provided by an observer) of a generalized Brunovsky state (\ie certain time derivatives of the flat output), which can be problematic in practice. In \cite{KolarRamsSchlacher:2017}, it is shown that under certain conditions, the measurements of a generalized Brunovsky state can be replaced by measurements of a (classical) state of the system. Based on our results, we will show that for $(x,u)$-flat systems with two inputs this is always possible.
		
		In Section \ref{se:not} we will introduce the notation used throughout this paper. In Section \ref{se:pre} we will state some properties of two-input $(x,u)$-flat systems which we will need in Section \ref{se:main}, where we present our main results. Section \ref{se:app} is dedicated to a practical application in flatness based tracking control of the rather theoretical main results.
	\section{Notation}\label{se:not}
		Let $\mathcal{X}$ be an $n$-dimensional smooth manifold, equipped with local coordinates $x^i$, $i=1,\ldots,n$. Its tangent bundle and cotangent bundle are denoted by $(\mathcal{T}(\mathcal{X}),\tau_\mathcal{X},\mathcal{X})$ and $(\mathcal{T}^\ast(\mathcal{X}),\tau^\ast_\mathcal{X},\mathcal{X})$. For these bundles we have the induced local coordinates $(x^i,\dot{x}^i)$ and $(x^i,\dot{x}_i)$ with respect to the bases $\{\partial_{x^i}\}$ and $\{\D x^i\}$, respectively. By $\partial_xf$ we denote the $m\times n$ Jacobian matrix of $f=(f^1,\ldots,f^m)$ with respect to $x=(x^1,\ldots,x^n)$ and by $\D\omega$ we denote the exterior derivative of a $p$-form $\omega$. The $k$-fold Lie-derivative of a function $\varphi$ along a vector field $v$ is denoted by $\Lie_v^k\varphi$. We make use of the Einstein summation convention. We write $\varphi^j_\alpha$ for the $\alpha$-th time derivative of $\varphi^j$ and use multi-indices to keep expressions involving many time derivatives short. Let $A=(a_1,\ldots,a_m)$ and $B=(b_1,\ldots,b_m)$ be two multi-indices with $a_j\leq b_j$, $j=1,\ldots,m$, which we abbreviate by $A\leq B$. Then
		\begin{align*}
			\begin{aligned}
				\varphi_{[A,B]}&=(\varphi^1_{[a_1,b_1]},\ldots,\varphi^m_{[a_m,b_m]})
			\end{aligned}
		\end{align*}
		where $\varphi^j_{[a_j,b_j]}=(\varphi^j_{a_j},\ldots,\varphi^j_{b_j})$ and $\varphi_{[0,A]}=\varphi_{[A]}$. Addition and subtraction of multi-indices is done component wise and $A\pm c=(a_1\pm c,\ldots,a_m\pm c)$ with an integer $c$, and $\#A=\sum_{j=1}^{m}a_j$. By $R=(r_1,r_2)$ we denote the unique multi-index associated to a flat output of a system with two inputs, where $r_j$ denotes the order of the highest derivative of $y^j$ needed to parameterize $x$ and $u$ by this flat output.
	\section{Preliminaries}\label{se:pre}
		In the following we work on a manifold $\mathcal{X}\times\mathcal{U}_{[l_u]}$ with coordinates $(x,u,u_1,\ldots,u_{l_u})$, where $u_\alpha$ denotes the $\alpha$-th time derivative of the input $u$ and $l_u$ is some large enough but finite integer. Consider a function $\varphi(x,u,u_1,\ldots,u_\gamma)$ on $\mathcal{X}\times\mathcal{U}_{[l_u]}$. As long as $\gamma<l_u$, \ie as long as $l_u$ is chosen big enough, the time derivative of this function is given by the Lie-derivative along the vector field
		\begin{align}\label{eq:2}
			\begin{aligned}
				f_u&=f^i(x,u)\partial_{x^i}+\sum_{\alpha=0}^{l_u-1}u^j_{\alpha+1}\partial_{u^j_\alpha}\,.
			\end{aligned}
		\end{align}
		In the following we assume that $l_u$ is chosen big enough such that $f_u$ acts as time derivative on all functions involved. 
		
		Consider a two-input $(x,u)$-flat system with a $(x,u)$-flat output $y=\varphi(x,u)$ and let us introduce a multi-index $K=(k_1,k_2)$, such that 
		\begin{align*}
			\begin{aligned}
				\Lie_{f_u}^{k_j-1}\varphi^j&=\varphi^j_{k_j-1}(x)\,,&&&\Lie_{f_u}^{k_j}\varphi^j&=\varphi^j_{k_j}(x,u)\,,
			\end{aligned}
		\end{align*}
		\ie $k_j$ denotes the relative degree of the component $\varphi^j$ of the flat output. Note that for a $(x,u)$-flat output where $\varphi^j=\varphi^j(x,u)$ actually depends on $u$, we have $k_j=0$. Furthermore, similar as in \cite{KolarSchoberlSchlacher:2015} and \cite{Kolar:2017}, let us introduce the codistributions
		\begin{align*}
			\begin{aligned}
				B_A&=\Span{\D \varphi_{[A]}}\cap\Span{\D x,\D u}\,.
			\end{aligned}
		\end{align*}
		For $A\leq K$, because of $\Span{\D \varphi_{[A]}}\subset\Span{\D x,\D u}$, we have $B_A=\Span{\D \varphi_{[A]}}$. Moreover, since $\Span{\D x,\D u}\subset\Span{\D \varphi_{[R]}}$, we have $B_{R}=\Span{\D x,\D u}$. These codistributions form the sequence
		\begin{align}\label{eq:3}
			\begin{aligned}
				B_K&\subset B_{K+1}\subset\ldots\subset B_R\,.
			\end{aligned}
		\end{align}
		(In Appendix \ref{app:a}, this sequence is illustrated for a simple model of a vehicle.) In the following, we exclude static feedback linearizable systems by assuming $R>K$. The following lemma provides a relation between the dimensions of the codistributions within the sequence \eqref{eq:3}.
		\begin{lem}\label{lem:1}
			For all $\beta\geq 0$ with $K+\beta<R$, we have
			\begin{align}\label{eq:4}
				\begin{aligned}
					\Dim{B_{K+\beta+1}}&=\Dim{B_{K+\beta}}+1\,.
				\end{aligned}
			\end{align}
		\end{lem}
		
		\begin{pf}
			Due to the functional independence of time derivatives of a flat output, the codistribution $B_{K+1}$ is given by
			\begin{align*}
				\begin{aligned}
					B_{K+1}&=(B_K\oplus\Span{\D \varphi^1_{k_1+1},\D \varphi^2_{k_2+1}})\cap\Span{\D x,\D u}\,.
				\end{aligned}
			\end{align*}
			The differentials $\D \varphi_{K+1}$ are of the form
			\begin{align}\label{eq:5}
				\begin{aligned}
					\D \varphi^j_{k_j+1}&=\ldots+\partial_{u^l}\varphi^j_{k_j}\D u^l_1\,,
				\end{aligned}
			\end{align}
			where the case that all four coefficients $\partial_{u^l}\varphi^j_{k_j}$ are zero can not occur due to the definition of $K$. Let us assume that $\Dim{B_{K+1}}=\Dim{B_K}$ holds, which means that there does not exist a linear combination of the differentials \eqref{eq:5} such that the differentials $\D u_1$ cancel out, \ie the $2\times 2$ Jacobian matrix $\partial_{u}\varphi_{K}$ is regular. Since the differentials of time derivatives of the functions $\varphi^j_{k_j}$ are of the form
			\begin{align*}
				\begin{aligned}
					\D\varphi^j_{k_j+\beta}&=\D\Lie _{f_u}^\beta\varphi^j_{k_j}=\ldots+\partial_{u^l}\varphi^j_{k_j}\D u^l_{\beta}\,,
				\end{aligned}
			\end{align*}
			with the same coefficients $\partial_{u^l}\varphi^j_{k_j}$ as in \eqref{eq:5}, also from the differentials $\varphi_{[K+\beta]}$, for arbitrary large $\beta$, no linear combinations contained in $\Span{\D x,\D u}$, which are not already contained in $B_K$, can be constructed. Thus, the assumption $\Dim{B_{K+1}}=\Dim{B_K}$ would imply $\Dim{B_{K+\beta}}=\Dim{B_K}$ for all $\beta\geq 0$, but this is in contradiction to $B_{R}=\Span{\D x,\D u}$. Thus, the $2\times 2$ Jacobian matrix $\partial_{u}\varphi_{K}$ must be singular. On the other hand, by the definition of $K$, none of the differentials $\D \varphi^j_{k_j+\beta}$, $\beta\geq 1$ on its own is contained in $\Span{\D x, \D u}$, therefore, as soon as we deal with codistributions with indices $K+\beta$, $\beta\geq 0$, $K+\beta<R$, the growth in dimension is exactly one with each increment of the index.\hfill$\square$
		\end{pf}
		Another consequence of the fact that none of the differentials $\D \varphi^j_{k_j+\beta}$, $\beta\geq 1$ on its own is contained in $\Span{\D x, \D u}$ is that the differences $r_1-k_1=r_2-k_2$ are always equal.
		The following lemma follows immediately form the proof of Lemma \ref{lem:1}.
		\begin{lem}\label{lem:2}
			For every $(x,u)$-flat output $y=\varphi(x,u)$ of the two-input system \eqref{eq:1}, $\Rank{\partial_u\varphi_K}=1$ holds.
		\end{lem}
		\begin{pf}
			In the proof of Lemma \ref{lem:1}, it turned out that the $2\times 2$ Jacobian matrix $\partial_u\varphi_K$ is singular and that it is not the zero matrix, thus its rank is one.\hfill$\square$
		\end{pf}
	\section{Main results}\label{se:main}
		In the following two sections, we will first show that every $(x,u)$-flat system with two inputs can be rendered static feedback linearizable by a special sub-class of endogenous dynamic feedback, namely prolongations of a single (new) input after a suitable static feedback has been applied. Then we will use this result to show that if a $(x,u)$-flat output of a two-input system is known, a so called generalized Brunovsky state for this system, with certain properties that are very useful for flatness based tracking control, can be constructed systematically. 
		\subsection{Linearization by prolongations}
			In the following, we make use of the static feedback
			\begin{align}\label{eq:6}
				\begin{aligned}
					\bar{u}^1&=\varphi^1_{k_1}(x,u)\\
					\bar{u}^2&=g(x,u)\,,
				\end{aligned}
			\end{align}
			\ie we define a derivative of the flat output as new input and choose $g(x,u)$ such that the Jacobian matrix of the right hand side of \eqref{eq:6} with respect to $u$ is regular. We will see that after applying this static feedback, the system can be rendered static feedback linearizable just by prolonging the new input $\bar{u}^1$, \ie preintegrating $\bar{u}^1$ suitably often. First, let us analyze the functions $\varphi_{[R]}$, \ie the flat output and its time derivatives as functions of $x$, $u$ and time derivatives of $u$ in the new coordinates given by \eqref{eq:6}\footnote{A static feedback or input transformation $\bar{u}=\Phi_u(x,u)$ entails the complete transformation
			\begin{align*}
				\begin{aligned}
					\bar{u}^j&=\Phi_u^j\\
					\bar{u}^j_1&=\Lie_{f_u}\Phi_u^j\\
					&\vdotswithin{=}\\
					\bar{u}^j_{l_u}&=\Lie_{f_u}^{l_u}\Phi_u^j
				\end{aligned}
			\end{align*}
			on $\mathcal{X}\times\mathcal{U}_{[l_u]}$. In these new coordinates the vector field $f_u$ is given by
			\begin{align*}
				\begin{aligned}
					f_u&=\underbrace{\bar{f}^i(x,\bar{u})\partial_{x^i}+\sum_{\alpha=0}^{l_u-1}\bar{u}^j_{\alpha+1}\partial_{\bar{u}^j_\alpha}}_{\bar{f}_u}+\ldots\partial_{\bar{u}^1_{l_u}}+\ldots\partial_{\bar{u}^2_{l_u}}\,,
				\end{aligned}
			\end{align*}
			where $\bar{f}=f\circ\hat{\Phi}_u$ with the inverse $\hat{\Phi}_u$ of $\Phi_u$. The in general non-zero components in the $\partial_{\bar{u}^j_{l_u}}$-directions do not bother us as long as $l_u$ is chosen big enough, \ie we can solely work with $\bar{f}_u$.}.
			\begin{thm}\label{thm:1}
				After applying the static feedback \eqref{eq:6}, the flat output and its time derivatives up to the order $R$ in the new coordinates are given by
				\begin{align}\label{eq:7}
					\begin{aligned}
					y_{[K-1]}&=\bar{\varphi}_{[K-1]}(x)\\
					y_K&=\begin{bmatrix}
					\bar{u}^1\\\bar{\varphi}^2_{k_2}(x,\bar{u}^1)
					\end{bmatrix}\\
					y_{K+1}&=\begin{bmatrix}
					\bar{u}^1_1\\\bar{\varphi}^2_{k_2+1}(x,\bar{u}^1,\bar{u}^1_1)
					\end{bmatrix}\\
					&\vdotswithin{=}\\
					y_{R-1}&=\begin{bmatrix}
					\bar{u}^1_{r_1-k_1-1}\\\bar{\varphi}^2_{k_2+r_1-k_1-1}(x,\bar{u}^1,\bar{u}^1_1,\ldots,\bar{u}^1_{r_1-k_1-1})
					\end{bmatrix}\\
					y_{R}
					&=\begin{bmatrix}
					\bar{u}^1_{r_1-k_1}\\\bar{\varphi}^2_{k_2+r_1-k_1}(x,\bar{u}^1,\bar{u}^2,\bar{u}^1_1,\ldots,\bar{u}^1_{r_1-k_1})
					\end{bmatrix}\,,
					\end{aligned}
				\end{align}
				where $\bar{u}^2$ only occurs in the last line in $\bar{\varphi}^2_{k_2+r_1-k_1}$. Furthermore, the map \eqref{eq:7} is a diffeomorphism.
			\end{thm}
			\begin{pf}
				After applying the input transformation \eqref{eq:6}, the flat output and its time derivatives up to order $K$ in the new coordinates are given by
				\begin{align*}
					\begin{aligned}
						y_{[K-1]}&=\bar{\varphi}_{[K-1]}(x)=\varphi_{[K-1]}(x)\,,\\
						y_K&=\begin{bmatrix}
						\bar{u}^1\\\bar{\varphi}^2_{k_2}(x,\bar{u}^1)
						\end{bmatrix}\,,
					\end{aligned}
				\end{align*}
				\ie the time derivatives up to order $K-1$ are not effected by this input transformation, since they do not depend on the inputs, and the time derivatives of order $K$, because of Lemma \ref{lem:2}, only depend on $\bar{u}^1$. Now let us assume that
				\begin{align*}
					\begin{aligned}
						y_{K+1}&=\begin{bmatrix}
						\bar{u}^1_1\\\bar{\varphi}^2_{k_2+1}(x,\bar{u}^1,\bar{u}^2,\bar{u}^1_1)
						\end{bmatrix}
					\end{aligned}
				\end{align*}
				and thus
				\begin{align*}
					\begin{aligned}
						y_{K+2}&=\begin{bmatrix}
						\bar{u}^1_2\\\bar{\varphi}^2_{k_2+2}(x,\bar{u}^1,\bar{u}^2,\bar{u}^1_1,\bar{u}^2_1,\bar{u}^1_2)
						\end{bmatrix}\,,
					\end{aligned}
				\end{align*}
				\ie we assume that $\bar{\varphi}^2_{k_2+1}$ actually depends on $\bar{u}^2$. Then in $\D\bar{\varphi}^2_{k_2+2}$, there necessarily occurs the differential $\D \bar{u}^2_1$. But this would imply that $\Dim{B_{K+2}}=\Dim{B_{K+1}}$ since there is no way to find a non-trivial linear combination of the differentials $\D\bar{\varphi}_{K+2}$ containing neither $\D \bar{u}^1_2$ nor $\D\bar{u}^2_1$. Thus, either $K+1=R$, \ie $B_{K+1}=\Span{\D x,\D u}$ is already the last codistribution of the sequence \eqref{eq:3}, or $\bar{\varphi}^2_{k_2+1}$ is actually independent of $\bar{u}^2$, \ie $\bar{\varphi}^2_{k_2+1}=\bar{\varphi}^2_{k_2+1}(x,\bar{u}^1,\bar{u}^1_1)$. Continuing this argumentation also for the functions $\bar{\varphi}^2_{k_2+\beta}$, $\beta\geq 2$, $k_2+\beta<r_2$, it follows that all of them are independent of $\bar{u}^2$, \ie that we actually have $\bar{\varphi}^2_{k_2+\beta}=\bar{\varphi}^2_{k_2+\beta}(x,\bar{u}^1,\bar{u}^1_1,\ldots,\bar{u}^1_{\beta})$, $\beta+k_2<r_2$ and $\bar{\varphi}^2_{r_2}=\bar{\varphi}^2_{r_2}(x,\bar{u}^1,\bar{u}^2,\bar{u}^1_1,\ldots,\bar{u}^1_{r_1-k_1})$ is the only function actually depending on $\bar{u}^2$. Thus, in conclusion (remember that $r_1-k_1=r_2-k_2$ holds and thus $k_2+r_1-k_1=r_2$) we have the form \eqref{eq:7}. 
				
				In the following we show that \eqref{eq:7} is actually a diffeomorphism $y_{[R]}=\hat{\bar{F}}_e(x,\bar{u}^1,\bar{u}^2,\bar{u}^1_{[1,r_1-k_1]})$ from the manifold $\mathcal{X}\times\mathcal{\bar{U}}_{[(0,r_1-k_1)]}$ with coordinates $(x,\bar{u}_{[(0,r_1-k_1)]})$ to the manifold $\mathcal{Y}_{[R]}$ with coordinates $y_{[R]}$. The functional independence of the right hand sides of \eqref{eq:7} follows directly from the fact that time derivatives of a flat output up to arbitrary order are functionally independent. What is left to show is that the number of variables on the right hand side coincides with the number of variables on the left hand side. From \eqref{eq:4}, it follows that $\Dim{B_R}=\Dim{B_K}+r_1-k_1$. With $\Dim{B_K}=\#K+2$ and $\Dim{B_R}=\Dim{\Span{\D x,\D u}}=n+2$ we thus have
				\begin{align*}
					\begin{aligned}
						n+2&=\#K+2+r_1-k_1
					\end{aligned}
				\end{align*}
				and from this, together with $r_1-k_1=r_2-k_2$, we obtain the two equations
				\begin{align}\label{eq:8}
					\begin{aligned}
						n-\#K&=r_1-k_1\\n-\#K&=r_2-k_2\,.
					\end{aligned}
				\end{align}
				Their sum yields 
				\begin{align*}
					\begin{aligned}
						2n-2\#K&=\#R-\#K&&\text{or}&n-\#K&=\#R-n\,.
					\end{aligned}
				\end{align*}
				Comparing the left hand sides of the latter equation and one of the equations in \eqref{eq:8}, it follows that $r_1-k_1=\#R-n$. Thus, for the map \eqref{eq:7}, we have $\Dim{x}+\Dim{u}+\Dim{\bar{u}^1_{[1,r_1-k_1]}}=n+2+\#R-n=\#R+2$ variables on the right hand side, which because of $\Dim{y_{[R]}}=\#R+2$ indeed coincides with the number of variables on the left hand side.\hfill$\square$
			\end{pf}
			The diffeomorphism \eqref{eq:7} is just the inverse of the parameterizing map of the prolonged system
			\begin{align}\label{eq:9}
				\begin{aligned}
					\dot{x}&=\bar{f}(x,\bar{u}^1,\bar{u}^2)\\
					\dot{\bar{u}}^1&=\bar{u}^1_1\\
					&\vdotswithin{=}\\
					\dot{\bar{u}}^1_{r_1-k_1-1}&=\bar{u}^1_{r_1-k_1}
				\end{aligned}
			\end{align}
			with the input $(\bar{u}^1_{r_1-k_1},\bar{u}^2)$, with respect to the flat output $y=\varphi(x,u)$, and since the parameterizing map of this prolonged system is a diffeomorphism, it is static feedback linearizable and $y=\varphi(x,u)$ is a linearizing output of it. These considerations can be summarized in the following corollary. 
			\begin{cor}\label{cor:1}
				Every $(x,u)$-flat system with two inputs can be rendered static feedback linearizable by $\#R-n$ prolongations of a suitably chosen (new) input (here $\bar{u}^1$).
			\end{cor}
		\subsection{Generalized Brunovsky states}
			For a detailed treatise on generalized states we refer to \cite{DelaleauRudolph:1998}, here, we will only outline the essentials. Based on the results of the previous section, in the following we will construct a generalized Brunovsky state for the $(x,u)$-flat two-input system \eqref{eq:1} that will fulfill additional properties, which we will need in Section \ref{se:app}. First, let us introduce the notion of a generalized state. As in \cite{KolarRamsSchlacher:2017}, we call an $n$-tuple $\tilde{x}=(\tilde{x}^1,\ldots,\tilde{x}^n)$ of functions
			\begin{align}\label{eq:10}
				\begin{aligned}
					\tilde{x}^i&=\tilde{\Phi}^i(x,u,u_1,\ldots,u_\gamma)
				\end{aligned}
			\end{align}
			with a regular Jacobian matrix $\partial_x\tilde{\Phi}$ a generalized state of \eqref{eq:1} and call the relation \eqref{eq:10} a generalized state transformation. A generalized Brunovsky state is a generalized state of the special form
			\begin{align*}
				\begin{aligned}
					\tilde{x}_B&=(y^1_{[\kappa_1-1]},\ldots,y^m_{[\kappa_m-1]})
				\end{aligned}
			\end{align*}
			with $\#\kappa=n$ where $\kappa=(\kappa_1,\ldots,\kappa_m)$ (the flat output and its time derivatives are indeed functions of $x$, $u$ and time derivatives of $u$). For flat systems such a generalized Brunovsky state always exists, see \cite{DelaleauRudolph:1998}. The following theorem states that for $(x,u)$-flat systems with two inputs, there always exists a generalized Brunovsky state which fulfills additional properties, which are useful for the design of tracking control laws as we will see in Section \ref{se:app}.			
			\begin{thm}\label{thm:2}
				The $n$-tuple
				\begin{align*}
					\begin{aligned}
						\tilde{x}_B&=y_{[\kappa-1]}
					\end{aligned}
				\end{align*}
				with $\kappa=(k_1,r_2)$ is a generalized Brunovsky state of the $(x,u)$-flat two-input system \eqref{eq:1}. For this generalized Brunovsky state, the properties $\kappa\leq R$ and $\Rank{\partial_{\tilde{x}_B}F_x(y_{[R-1]})}=n$, hold.
			\end{thm}
			\begin{pf}
				To show that the $n$ components
				\begin{align*}
					\begin{aligned}
						y_{[\kappa-1]}&=(\underbrace{y,y_1,\ldots,y_{K-1}}_{\#K},\underbrace{y^2_{k_2},y^2_{k_2+1},\ldots,y^2_{r_2-1}}_{r_2-k_2})
					\end{aligned}
				\end{align*}
				of $y_{[R-1]}$, $\#\kappa=\#K+r_2-k_2=n$ follows from \eqref{eq:8}, form a generalized Brunovsky state of the system, we have to show that $\tilde{x}_B=\bar{\varphi}_{[\kappa-1]}$ is a generalized state transformation \eqref{eq:10}, \ie that the Jacobian matrix $\partial_x\bar{\varphi}_{[\kappa-1]}$ is regular. This proof, as well as the proof of the additional property $\Rank{\partial_{\tilde{x}_B}F_x(y_{[R-1]})}=n$, is based on the diffeomorphism \eqref{eq:7}. This diffeomorphism contains the diffeomorphism $y_{[R-1]}=\hat{\bar{F}}_{e,red}(x,\bar{u}^1_{[r_1-k_1-1]})$ from the manifold \mbox{$\mathcal{X}\times\mathcal{\bar{U}}^1_{[r_1-k_1-1]}$} with coordinates $(x,\bar{u}^1_{[r_1-k_1-1]})$ to the manifold $\mathcal{Y}_{[R-1]}$ with coordinates $y_{[R-1]}$, and is given by
				\begin{align}\label{eq:11}
					\begin{aligned}
						y_{[K-1]}&=\bar{\varphi}_{[K-1]}(x)\\
						y_K&=\begin{bmatrix}
						\bar{u}^1\\\bar{\varphi}^2_{k_2}(x,\bar{u}^1)
						\end{bmatrix}\\
						y_{K+1}&=\begin{bmatrix}
						\bar{u}^1_1\\\bar{\varphi}^2_{k_2+1}(x,\bar{u}^1,\bar{u}^1_1)
						\end{bmatrix}\\
						&\vdotswithin{=}\\
						y_{R-1}&=\begin{bmatrix}
						\bar{u}^1_{r_1-k_1-1}\\\bar{\varphi}^2_{k_2+r_1-k_1-1}(x,\bar{u}^1,\bar{u}^1_1,\ldots,\bar{u}^1_{r_1-k_1-1})
						\end{bmatrix}\,.
					\end{aligned}
				\end{align}
				Its inverse can be considered as the parameterizing map for the state $(x,\bar{u}^1_{[r_1-k_1-1]})$ of the prolonged system \eqref{eq:9} and is given by
				\begin{align}\label{eq:12}
					\begin{aligned}
						x&=F_x(y_{[R-1]})\\
						\bar{u}^1&=y^1_{k_1}\\
						\bar{u}^1_1&=y^1_{k_1+1}\\
						&\vdotswithin{=}\\
						\bar{u}^1_{r_1-k_1-1}&=y^1_{r_1-1}\,.
					\end{aligned}
				\end{align}
				The proof of the regularity of $\partial_x\bar{\varphi}_{[\kappa-1]}$ is based on the Jacobian matrix of \eqref{eq:11} and can be found in Appendix \ref{app:b}. It is similar to the proof of the property $\Rank{\partial_{\tilde{x}_B}F_x(y_{[R-1]})}=n$, which is based on the Jacobian matrix of \eqref{eq:12}, and which we present here. The Jacobian matrix of \eqref{eq:12} reads
				\begin{align}\label{eq:13}
					\begin{aligned}
						&\partial_{y_{[R-1]}}\bar{F}_{e,red}=\\
						&\begin{bmatrix}
						\partial_yF_x&\partial_{y_1}F_x&\ldots&\partial_{y_{K-1}}F_x&\partial_{y_{K}}F_x&\partial_{y_{K+1}}F_x&\ldots&\partial_{y_{R-1}}F_x\\
						0&0&\ldots&0&\begin{bmatrix}1&0\end{bmatrix}&0&\ldots&0\\
						0&0&\ldots&0&0&\begin{bmatrix}1&0\end{bmatrix}&\ldots&0\\
						\vdots&\vdots&&\vdots&\vdots&\vdots&&\vdots\\
						0&0&\ldots&0&0&0&\ldots&\begin{bmatrix}1&0\end{bmatrix}
						\end{bmatrix}
					\end{aligned}
				\end{align}
				and since it is the Jacobian matrix of a diffeomorphism, its columns (as well as its rows) are linearly independent. For the Jacobian matrix $\partial_{y_{[\kappa-1]}}\bar{F}_{e,red}$ we obtain
				\begin{align}\label{eq:14}
					\begin{aligned}
					\partial_{y_{[\kappa-1]}}\bar{F}_{e,red}&=\begin{bmatrix}
							\partial_{y_{[\kappa-1]}}F_x\\
							0
						\end{bmatrix}
					\end{aligned}
				\end{align}
				and the columns of this matrix are just certain columns of \eqref{eq:13}, which in turn are all linearly independent from each other. The block of zeros underneath the Jacobian matrix $\partial_{y_{[\kappa-1]}}F_x$ in \eqref{eq:14} does not contribute to the rank of $\partial_{y_{[\kappa-1]}}\bar{F}_{e,red}$, thus the $n$ columns of the $n\times n$ Jacobian matrix $\partial_{y_{[\kappa-1]}}F_x$ are linearly independent, which means $\partial_{y_{[\kappa-1]}}F_x$ is regular.\hfill$\square$
			\end{pf}
			(The rather theoretical results of this section are also illustrated on the vehicle model in Appendix \ref{app:a}.)
			\begin{rem}\label{rem:1}
				In \eqref{eq:6}, instead of choosing $\bar{u}^1=\varphi^1_{k_1}(x,u)$ one could also choose $\bar{u}^1=\varphi^2_{k_2}(x,u)$ and proceeding with this choice, it follows that $\tilde{x}_B=y_{[\kappa-1]}$ with $\kappa=(r_1,k_2)$ is also a valid generalized Brunovsky state which also possesses the properties as in Theorem \ref{thm:2}.
			\end{rem}
		
	\section{Application}\label{se:app}
		The above results are useful for flatness based tracking control. In \cite{DelaleauRudolph:1998} a method for tracking control based on exact linearization by a quasi-static state feedback is presented. This approach, in contrast to tracking control based on exact linearization by an endogenous dynamic feedback, yields a static control law, but it requires measurements of a generalized Brunovsky state (\ie measurements of the flat output and time derivatives of the flat output up to a certain order). In \cite{KolarRamsSchlacher:2017} a method which allows the use of measurements of the state of the system instead of measurements of the generalized Brunovsky state is presented. An open question is whether this method is always applicable. In the following we will briefly cover tracking control based on quasi-static state feedback, for details see \cite{DelaleauRudolph:1998}. Then we will recapitulate the method presented in \cite{KolarRamsSchlacher:2017} and will see that for $(x,u)$-flat systems with two inputs it is always applicable.
		
		Flat systems can be exactly linearized by a quasi-static state feedback, see \cite{DelaleauRudolph:1998}. If a generalized Brunovsky state is measured, the system can be exactly linearized by the quasi-static state feedback
		\begin{align}\label{eq:15}
			\begin{aligned}
				u&=F_u(\tilde{x}_B,v_{[R-\kappa]})\,,
			\end{aligned}
		\end{align}
		which is constructed from the map $F_u$ by replacing $y_{[\kappa-1]}$ by $\tilde{x}_B$ and $y_{[\kappa,R]}$ by the new input $v=y_\kappa$ and its time derivatives up to order $R-\kappa$. The feedback
		\begin{align}\label{eq:16}
			\begin{aligned}
				v^j&=y^{j,d}_{\kappa_j}-\sum_{\beta=0}^{\kappa_j-1}a^{j,\beta}\left(y^j_\beta-y^{j,d}_\beta\right)\,,
			\end{aligned}
		\end{align}
		where $y^d$ denotes the reference trajectory, results in the linear tracking error dynamics
		\begin{align}\label{eq:17}
			\begin{aligned}
				e^j_{\kappa_j}+\sum_{\beta=0}^{\kappa_j-1}a^{j,\beta}e^j_\beta&=0\,,
			\end{aligned}
		\end{align}
		where $e^j=y^j-y^{j,d}$, and the roots of the characteristic polynomials can be adjusted by the coefficients $a^{j,\beta}\in\mathbb{R}$. The control law \eqref{eq:15} contains the time derivatives $v_{[R-\kappa]}$ of $v$, but those can be eliminated. The $\lambda$-th time derivative $v_\lambda$ of $v$ is given by 
		\begin{align*}
			\begin{aligned}
				v^j_\lambda&=y^{j,d}_{\kappa_j+\lambda}-\sum_{\beta=0}^{\kappa_j-1}a^{j,\beta}\left(y^j_{\beta+\lambda}-y^{j,d}_{\beta+\lambda}\right),\,\lambda=1,\ldots,r_j-\kappa_j,
			\end{aligned}
		\end{align*}
		which follows directly from differentiating \eqref{eq:16}. Now replacing the time derivatives $y_{[\kappa,R-1]}$ of the flat output which are not contained in $\tilde{x}_B$ by $v_{[R-\kappa-1]}$ results in the system of linear equations
		\begin{align}\label{eq:18}
			\begin{aligned}
				v^j&=y^{j,d}_{\kappa_j}-\sum_{\beta=0}^{\kappa_j-1}a^{j,\beta}\left(y^j_\beta-y^{j,d}_\beta\right)\\
				v^j_1&=y^{j,d}_{\kappa_j+1}-a^{j,\kappa_j-1}\left(v^j-y^{j,d}_{\kappa_j}\right)-\\&\hspace{12ex}\sum_{\beta=0}^{\kappa_j-2}a^{j,\beta}\left(y^j_{\beta+1}-y^{j,d}_{\beta+1}\right)\\
				&\vdotswithin{=}\\
				v^j_{r_j-\kappa_j}&=\ldots\,,
			\end{aligned}
		\end{align}
		which can be solved systematically from top to bottom for the unknowns $v_{[R-\kappa]}$. Substituting the solution for $v_{[R-\kappa]}$ into \eqref{eq:15} gives a control law of the form
		\begin{align*}
			\begin{aligned}
				u&=\alpha(\tilde{x}_B,y^d_{[R]})\,.
			\end{aligned}
		\end{align*}
		This kind of tracking control requires measurements of the generalized Brunovsky state $\tilde{x}_B$. The main idea of the method presented in \cite{KolarRamsSchlacher:2017} to obtain a control law independent of $\tilde{x}_B$, is to replace $\tilde{x}_B$ by solving $x=F_x(y_{[R-1]})$ for $\tilde{x}_B$. According to Theorem \ref{thm:2}, $\tilde{x}_B=y_{[\kappa-1]}$ with $\kappa=(k_1,r_2)$ is a generalized Brunovsky state of the two-input system \eqref{eq:1} with $(x,u)$-flat output $y=\varphi(x,u)$ and for this generalized Brunovsky state the Jacobian matrix $\partial_{\tilde{x}_B}F_x(y_{[R-1]})$ is regular. The implicit function theorem then guarantees that $x=F_x(y_{[R-1]})$ can indeed locally be solved for $\tilde{x}_B$ as function of $x$ and $y_{[\kappa,R-1]}$, \ie we obtain a relation of the form
		\begin{align}\label{eq:19}
			\begin{aligned}
				\tilde{x}_B&=\phi(x,y_{[\kappa,R-1]})
			\end{aligned}
		\end{align} 
		and by replacing $y_{[\kappa,R-1]}$ in \eqref{eq:19} by $v_{[R-\kappa-1]}$ we obtain
		\begin{align}\label{eq:20}
			\begin{aligned}
				\tilde{x}_B&=\tilde{\phi}(x,v_{[R-\kappa-1]})\,.
			\end{aligned}
		\end{align}
		The relation \eqref{eq:20} is of the form
		\begin{align}\label{eq:21}
			\begin{aligned}
				y^1_{[\kappa_1-1]}&=\bar{\varphi}^1_{[k_1-1]}(x)\\
				y^2_{[\kappa_2-1]}&=(\bar{\varphi}^2_{[k_2-1]}(x),\bar{\varphi}^2_{[k_2,r_2-1]}(x,v^1_{[r_2-k_2-1]}))\,,
			\end{aligned}
		\end{align}	
		\ie the components $y^1_{[\kappa_1-1]}$ only depend on $x$, and $v^2$ or time derivatives of $v^2$ do not occur at all. These properties follow from $\kappa_1-1<k_1$ and from $r_2-\kappa_2-1<0$.
		
		Substituting the relation \eqref{eq:20} into \eqref{eq:15} gives a control law of the form
		\begin{align}\label{eq:22}
			\begin{aligned}
			u&=F_u(\tilde{\phi}(x,v_{[R-\kappa-1]}),v_{[R-\kappa]})\,.
			\end{aligned}
		\end{align}
		As before, we want to eliminate the time derivatives $v_{[R-\kappa]}$ from the control law. For that, we replace $y_{[\kappa-1]}$ in \eqref{eq:18} by \eqref{eq:20}, respecting the special structure \eqref{eq:21} of \eqref{eq:20}, to obtain the equation system
		\begin{align}\label{eq:23}
			\begin{aligned}
				v^1&=y^{1,d}_{k_1}-\sum_{\beta=0}^{k_1-1}a^{1,\beta}\left(\bar{\varphi}^1_\beta(x)-y^{1,d}_\beta\right)\\
				v^1_1&=y^{1,d}_{k_1+1}-a^{1,k_1-1}\left(v^1-y^{1,d}_{k_1}\right)-\\
				&\hspace{12ex}\sum_{\beta=0}^{k_1-2}a^{1,\beta}\left(\bar{\varphi}^1_{\beta+1}(x)-y^{1,d}_{\beta+1}\right)\\
				&\vdotswithin{=}\\
				v^1_{r_1-k_1}&=\ldots\\
				v^2&=y^{2,d}_{r_2}-\sum_{\beta=0}^{r_2-1}a^{2,\beta}\left(\bar{\varphi}^2_\beta(x,v^1,\ldots,v^1_{\beta-k_2})-y^{2,d}_\beta\right)\,,
			\end{aligned}
		\end{align}
		which again can be solved systematically from top to bottom for the unknowns $v_{[R-\kappa]}$ as function of $x$ and $y^d_{[R]}$. Inserting this solution for $v_{[R-\kappa]}$ into \eqref{eq:22} yields a control law of the desired form
		\begin{align*}
			\begin{aligned}
				u&=\alpha(x,y^d_{[R]})\,,
			\end{aligned}
		\end{align*}
		\ie a control law which only depends on $x$ and $y^d_{[R]}$ and results in the linear tracking error dynamics \eqref{eq:17}.
		
		Roughly speaking, the above control law follows from the parameterizing map $u=F_u(y_{[R]})$ by replacing $y_{[\kappa-1]}$ by $\tilde{x}_B$, which in turn is expressed as function of $x$ and $v_{[R-\kappa-1]}$, and replacing $y_{[\kappa,R]}$ by $v_{[R-\kappa]}$. The time derivatives $v_{[R-\kappa]}$ are then expressed as functions of $x$ and $y^d_{[R]}$, \ie the solution of the equation system \eqref{eq:23}. In conclusion, we replace $y_{[R]}$ in $u=F_u(y_{[R]})$ by $x$ and $y^d_{[R]}$. In Appendix \ref{app:a}, we derive such a control law for the planar VTOL aircraft, which is also treated \eg in \cite{FliessLevineMartinRouchon:1999}, \cite{SchoberlRiegerSchlacher:2010} or \cite{SchoberlSchlacher:2011}.
		\begin{rem}\label{rem:2}
			Instead of actually solving the parameterization $x=F_x(y_{[R-1]})$ for $\tilde{x}_B$ as function of $x$ and $y_{[\kappa,R-1]}$ and then replacing $y_{[\kappa,R-1]}$ by $v_{[R-\kappa-1]}$, there is another way to obtain \eqref{eq:20}. The generalized Brunovsky state $\tilde{x}_B=y_{[\kappa-1]}$ is just a selection of certain time derivatives of the flat output. Because of $\kappa_1=k_1$, the components $y^1_{[\kappa_1-1]}=\varphi^1_{[\kappa_1-1]}(x)$ are readily obtained as function of $x$ just by computing the corresponding time derivatives up to the order $\kappa_1-1$. The same holds for the components $y^2_{[k_2-1]}=\varphi^2_{[k_2-1]}(x)$. By further time differentiation we obtain the expressions for the remaining components $y^2_{[k_2,r_2-1]}$, but in those, besides $x$, also the inputs and time derivatives of the inputs occur. However, if the input transformation \eqref{eq:6} is applied, because of \eqref{eq:7}, only $\bar{u}^1$ and time derivatives of $\bar{u}^1$ occur. Also from \eqref{eq:7}, we can read off $\bar{u}^1=y^1_{k_1}$ and since $\kappa_1=k_1$ and $y^1_{\kappa_1}=v^1$, we have $\bar{u}^1=v^1$. Thus, all we have to do to obtain the components $y^2_{[k_2,r_2-1]}$ as functions of $x$ and $v_{[R-\kappa-1]}$, is to replace $\bar{u}^1$ and time derivatives of $\bar{u}^1$ by $v^1$ and time derivatives of $v^1$, \ie insert $\bar{u}^1_\alpha=v^1_\alpha$ for $\alpha=0,\ldots,r_2-k_2-1$.
		\end{rem}
		\begin{rem}\label{rem:3}
			The open questions concerning the applicability of the method presented in \cite{KolarRamsSchlacher:2017} for the general case of flat systems with $m>2$ inputs are the existence of a generalized Brunovsky state for which there exists a relation of the form \eqref{eq:20}, and, if such a relation exists, the solvability of the system of non-linear equations obtained by inserting this relation into \eqref{eq:18}\,.
		\end{rem}
		\vspace{-0.9ex}
		\section{Conclusion}
		\vspace*{-2ex}
			We have shown that $(x,u)$-flat systems with two inputs can be rendered static feedback linearizable by a special sub-class of endogenous dynamic feedback. Future research will cover the question whether also for $(x,u)$-flat systems with $m>2$ inputs a restriction to certain classes of endogenous dynamic feedbacks is possible. Furthermore, we have presented an application of this theoretical result to the design of flatness based tracking controls. We have shown that the method presented in \cite{KolarRamsSchlacher:2017} for the design of a quasi-static tracking control, that uses measurements of a classical state instead of measurements of a generalized Brunovsky state, is always applicable. Future research will cover the question whether this method is also generally applicable for $(x,u)$-flat systems with $m>2$ inputs.
			\vspace{-0.6ex}
		\bibliography{Bibliography}           

\begin{thebibliography}{19}
\providecommand{\natexlab}[1]{#1}
\providecommand{\url}[1]{\texttt{#1}}
\providecommand{\urlprefix}{URL }
\expandafter\ifx\csname urlstyle\endcsname\relax
  \providecommand{\doi}[1]{doi:\discretionary{}{}{}#1}\else
  \providecommand{\doi}{doi:\discretionary{}{}{}\begingroup
  \urlstyle{rm}\Url}\fi

\bibitem[{Charlet and L{\'e}vine(1989)}]{CharletLevine:1989}
Charlet, B. and L{\'e}vine, J. (1989).
\newblock On dynamic feedback linearization.
\newblock \emph{Systems $\&$ Control Letters}, 13, 143--151.

\bibitem[{Charlet et~al.(1991)Charlet, L{\'e}vine, and
  Marino}]{CharletLevineMarino:1991}
Charlet, B., L{\'e}vine, J., and Marino, R. (1991).
\newblock Sufficient conditions for dynamic state feedback linearization.
\newblock \emph{SIAM J. Control Optim.}, 29(1), 38--57.

\bibitem[{Delaleau and Rudolph(1998)}]{DelaleauRudolph:1998}
Delaleau, E. and Rudolph, J. (1998).
\newblock Control of flat systems by quasi-static feedback of generalized
  states.
\newblock \emph{International Journal of Control}, 71(5), 745--765.

\bibitem[{Fliess et~al.(1992)Fliess, L{\'e}vine, Martin, and
  Rouchon}]{FliessLevineMartinRouchon:1992}
Fliess, M., L{\'e}vine, J., Martin, P., and Rouchon, P. (1992).
\newblock Sur les syst{\`e}mes non lin{\'e}aires diff{\'e}rentiellement plats.
\newblock \emph{Comptes rendus de l'Acad{\'e}mie des sciences. S{\'e}rie I,
  Math{\'e}matique}, 315, 619--624.

\bibitem[{Fliess et~al.(1995)Fliess, L{\'e}vine, Martin, and
  Rouchon}]{FliessLevineMartinRouchon:1995}
Fliess, M., L{\'e}vine, J., Martin, P., and Rouchon, P. (1995).
\newblock Flatness and defect of non-linear systems: introductory theory and
  examples.
\newblock \emph{International Journal of Control}, 61(6), 1327--1361.

\bibitem[{Fliess et~al.(1999)Fliess, L{\'e}vine, Martin, and
  Rouchon}]{FliessLevineMartinRouchon:1999}
Fliess, M., L{\'e}vine, J., Martin, P., and Rouchon, P. (1999).
\newblock A {L}ie-{B}{\"a}cklund approach to equivalence and flatness of
  nonlinear systems.
\newblock \emph{IEEE Transactions on Automatic Control}, 44(5), 922--937.

\bibitem[{Fossas and Franch(2000)}]{FossasFranch:2000}
Fossas, E. and Franch, J. (2000).
\newblock Linearization by prolongations and quotient of modules, a case study:
  the vertical take off and landing aircraft.
\newblock In \emph{Proceedings WSES International Conference on Applied and
  Theoretical Mathematics}, 1961--1965.

\bibitem[{Franch and Fossas(2005)}]{FranchFossas:2005}
Franch, J. and Fossas, E. (2005).
\newblock Linearization by prolongations: New bounds on the number of
  integrators.
\newblock \emph{European Journal of Control}, 11(2), 171 -- 179.

\bibitem[{Jakubczyk and Respondek(1980)}]{JakubczykRespondek:1980}
Jakubczyk, B. and Respondek, W. (1980).
\newblock On linearization of control systems.
\newblock \emph{Bull. Acad. Polonaise Sci. Ser. Sci. Math.}, 28, 517--522.

\bibitem[{Kolar(2017)}]{Kolar:2017}
Kolar, B. (2017).
\newblock \emph{Contributions to the Differential Geometric Analysis and
  Control of Flat Systems}.
\newblock Shaker Verlag, Aachen.

\bibitem[{Kolar et~al.(2015)Kolar, Sch{\"o}berl, and
  Schlacher}]{KolarSchoberlSchlacher:2015}
Kolar, B., Sch{\"o}berl, M., and Schlacher, K. (2015).
\newblock Remarks on a triangular form for 1-flat {P}faffian systems with two
  inputs.
\newblock In \emph{Proceedings 1st IFAC Conference on Modelling, Identification
  and Control of Nonlinear Systems (MICNON)}.
\newblock IFAC-PapersOnLine, volume 48, issue 11, pages 109--114.

\bibitem[{Kolar et~al.(2017)Kolar, Rams, and
  Schlacher}]{KolarRamsSchlacher:2017}
Kolar, B., Rams, H., and Schlacher, K. (2017).
\newblock Time-optimal flatness based control of a gantry crane.
\newblock \emph{Control Engineering Practice}, 60, 18--27.

\bibitem[{Nicolau and Respondek(2017)}]{NicolauRespondek:2017}
Nicolau, F. and Respondek, W. (2017).
\newblock Flatness of multi-input control-affine systems linearizable via
  one-fold prolongation.
\newblock \emph{SIAM J. Control and Optimization}, 55, 3171--3203.

\bibitem[{Nijmeijer and van~der Schaft(1990)}]{NijmeijervanderSchaft:1990}
Nijmeijer, H. and van~der Schaft, A. (1990).
\newblock \emph{Nonlinear Dynamical Control Systems}.
\newblock Springer, New York.

\bibitem[{Schlacher and Sch{\"o}berl(2013)}]{SchlacherSchoberl:2013}
Schlacher, K. and Sch{\"o}berl, M. (2013).
\newblock A jet space approach to check {P}faffian systems for flatness.
\newblock In \emph{Proceedings 52nd IEEE Conference on Decision and Control
  (CDC)}, 2576--2581.

\bibitem[{Sch{\"o}berl et~al.(2010)Sch{\"o}berl, Rieger, and
  Schlacher}]{SchoberlRiegerSchlacher:2010}
Sch{\"o}berl, M., Rieger, K., and Schlacher, K. (2010).
\newblock System parametrization using affine derivative systems.
\newblock In \emph{Proceedings 19th International Symposium on Mathematical
  Theory of Networks and Systems (MTNS)}, 1737--1743.

\bibitem[{Sch{\"o}berl and Schlacher(2011)}]{SchoberlSchlacher:2011}
Sch{\"o}berl, M. and Schlacher, K. (2011).
\newblock On calculating flat outputs for pfaffian systems by a reduction
  procedure - demonstrated by means of the vtol example.
\newblock In \emph{9th IEEE International Conference on Control $\&$ Automation
  (ICCA11)}, 477--482.

\bibitem[{Sch{\"o}berl and Schlacher(2014)}]{SchoberlSchlacher:2014}
Sch{\"o}berl, M. and Schlacher, K. (2014).
\newblock On an implicit triangular decomposition of nonlinear control systems
  that are 1-flat - a constructive approach.
\newblock \emph{Automatica}, 50, 1649--1655.

\bibitem[{Sluis and Tilbury(1996)}]{SluisTilbury:1996}
Sluis, W. and Tilbury, D. (1996).
\newblock A bound on the number of integrators needed to linearize a control
  system.
\newblock \emph{Systems \& Control Letters}, 29, 43--50.

\end{thebibliography}
		                                                  
	\clearpage
	\newpage
	\appendix
	\section{Examples}\label{app:a}
		In the following we will illustrate some of the theoretical results of the previous sections based on a simple model of a vehicle, taken from \cite{NijmeijervanderSchaft:1990}, and an academic example taken from \cite{SchoberlRiegerSchlacher:2010}. Furthermore, we will derive a tracking control law according to Section \ref{se:app} for the planar VTOL aircraft. In the following we use the abbreviations $\Cos(\cdot)=\cos(\cdot)$, $\Sin(\cdot)=\sin(\cdot)$ and $\Tan(\cdot)=\tan(\cdot)$.
		
			\subsubsection{Example 1.}\label{ex:1}The simple vehicle model is given by
			\begin{align*}
				\begin{aligned}
					\dot{x}^1&=\Sin(x^3)u^1\\
					\dot{x}^2&=\Cos(x^3)u^1\\
					\dot{x}^3&=u^2\,.
				\end{aligned}
			\end{align*}
			For this system, the vector field \eqref{eq:2} reads
			\begin{align*}
				\begin{aligned}
					f_u&=\Sin(x^3)u^1\partial_{x^1}+\Cos(x^3)u^1\partial_{x^2}+u^2\partial_{x^3}+\\
					&\hspace{8ex}\sum_{\alpha=0}^{l_u-1}(u^1_{\alpha+1}\partial_{u^1_\alpha}+u^2_{\alpha+1}\partial_{u^2_\alpha})\,.
				\end{aligned}
			\end{align*}
			The position $(x^1,x^2)$ is a $x$-flat output of this system, but for demonstration purposes, let us proceed with the $(x,u)$-flat output
			\begin{align*}
				\begin{aligned}
					y&=(\underbrace{x^1+\Cos(x^3)u^1}_{\varphi^1(x,u)},\underbrace{x^2}_{\varphi^2(x)})\,,
				\end{aligned}
			\end{align*}
			which is obtained from the $x$-flat output $\bar{y}=(x^1,x^2)$ by adding the time derivative of the second component to the first one \ie $y^1=\bar{y}^1+\bar{y}^2_1$, $y^2=\bar{y}^2$. One verifies that the parameterizing map with respect to this flat output is of the form
			\begin{align*}
				\begin{aligned}
					x&=F_x(y^1_{[1]},y^2_{[2]})\\
					u&=F_u(y^1_{[2]},y^2_{[3]})\,,
				\end{aligned}
			\end{align*}
			\ie $R=(2,3)$. The function $\varphi^1$ explicitly depends on an input, thus $k_1=0$. The function $\varphi^2$ is independent of the inputs, but
			\begin{align*}
				\begin{aligned}
					\varphi^2_1&=\Lie_{f_u}\varphi^2=\Cos(x^3)u^1
				\end{aligned}
			\end{align*}
			explicitly depends on an input, thus $k_2=1$ and $K=(0,1)$. The differentials $\D\varphi^1$, $\D\varphi^2$, $\D\varphi^2_1$ form a basis for the codistribution $B_K$ in \eqref{eq:3}, the complete sequence \eqref{eq:3} consists of the codistributions
			\begin{align*}
				\begin{aligned}
					B_K&=\Span{\D x^1,\D x^2,\Tan(x^3)u^1\D x^3-\D u^1}\\
					B_{K+1}&=\Span{\D x^1,\D x^2,\D x^3,\D u^1}\\
					B_{K+2}&=B_R=\Span{\D x,\D u}\,.
				\end{aligned}
			\end{align*}
			The input transformation \eqref{eq:6} can be chosen as
			\begin{align}\label{eq:24}
				\begin{aligned}
					\bar{u}^1&=\varphi^1_{k_1}=\varphi^1=x^1+\Cos(x^3)u^1\\
					\bar{u}^2&=u^2\,.
				\end{aligned}
			\end{align}
			In these coordinates the flat output and its time derivatives up to the order $R$ read
			\begin{align*}
				\begin{aligned}
					\bar{\varphi}_{[K-1]}:&&y^2&=x^2\\
					\bar{\varphi}_{K}:&&y^1&=\bar{u}^1\,,&&y^2_1=\bar{u}^1-x^1\\
					\bar{\varphi}_{K+1}:&&y^1_1&=\bar{u}^1_1\,,&&y^2_2=(x^1-\bar{u}^1)\Tan(x^3)+\bar{u}^1_1\\
					\bar{\varphi}_{K+2}=\bar{\varphi}_R:&&y^1_2&=\bar{u}^1_2\,,&&
					y^2_3=(\bar{u}^1-x^1)((1-\bar{u}^2)\Tan^2(x^3)\\
					&&&&&\hspace{9ex}-\bar{u}^2)-\bar{u}^1_1\Tan(x^3)+\bar{u}^1_2\,,
				\end{aligned}
			\end{align*}
			which is of the form \eqref{eq:7}, indeed, $\bar{u}^2$ only occurs in the last line in $y^2_3$. According to Theorem \ref{thm:2}, $\tilde{x}_B=y_{[\kappa-1]}$ with $\kappa=(k_1,r_2)=(0,3)$, \ie $\tilde{x}_B=y^2_{[2]}=(y^2,y^2_1,y^2_2)$ is a generalized Brunovsky state of the system. Indeed, the Jacobian matrix
			\begin{align*}
				\begin{aligned}
					\partial_x\bar{\varphi}^2_{[2]}&=\begin{bmatrix}
						0&1&0\\
						-1&0&0\\
						\Tan(x^3)&0&(x^1-\bar{u}^1)(1+\Tan^2(x^3))
					\end{bmatrix}
				\end{aligned}
			\end{align*}
			is regular. According to Theorem \ref{thm:2}, for this generalized Brunovsky state also the Jacobian matrix $\partial_{\tilde{x}_B}F_x$ is regular. One verifies that it is given by
			\begin{align*}
				\begin{aligned}
					\partial_{\tilde{x}_B}F_x&=\begin{bmatrix}
						~0&-1&0\\
						~1&0&0\\
						~0&\frac{y^2_2-y^1_1}{(y^1_1-y^2_2)^2+(y^2_1)^2}&\frac{-y^2_1}{(y^1_1-y^2_2)^2+(y^2_1)^2}
					\end{bmatrix}
				\end{aligned}
			\end{align*}
			and it is indeed regular.
			
			Note that according to Remark \ref{rem:1}, instead of the choice \eqref{eq:24}, one could also choose
				\begin{align*}
					\begin{aligned}
						\bar{u}^1=\varphi^2_{k_2}=\varphi^2_1=\Cos(x^3)u^1\,,
					\end{aligned}
				\end{align*}
				together with a suitable complement, \eg again $\bar{u}^2=u^2$, as input transformation. Furthermore, $\tilde{x}_B=y_{[\kappa-1]}$ with $\kappa=(r_1,k_2)=(2,1)$, \ie $\tilde{x}_B=(y^1_{[1]},y^2)=(y^1,y^1_1,y^2)$ is a generalized Brunovsky state of the system, which also fulfills that the Jacobian matrix $\partial_{\tilde{x}_B}F_x$ is regular.

			\subsubsection{Example 2.}\label{ex:2}
			Consider the system
			\begin{align*}
				\begin{aligned}
					\dot{x}^1&=u^1\\
					\dot{x}^2&=u^2\\
					\dot{x}^3&=\sqrt{u^1u^2}\,,
				\end{aligned}
			\end{align*}
			with the $(x,u)$-flat output
			\begin{align*}
				\begin{aligned}
					y&=(x^2-x^1u^2/u^1,x^3-x^1\sqrt{u^2/u^1})
				\end{aligned}
			\end{align*}
			(unlike the vehicle model, this system does not possess a $x$-flat output). The parameterizing map with respect to this flat output is of the form
			\begin{align*}
				\begin{aligned}
					x&=F_x(y^1_{[2]},y^2_{[2]})\\
					u&=F_u(y^1_{[3]},y^2_{[3]})\,,
				\end{aligned}
			\end{align*}
			\ie $R=(3,3)$ and since both components of the flat output $y$ explicitly depend on the inputs, we have $K=(0,0)$. The input transformation \eqref{eq:6} can be chosen as
			\begin{align*}
				\begin{aligned}
					\bar{u}^1&=x^2-x^1u^2/u^1\\
					\bar{u}^2&=u^2\,.
				\end{aligned}
			\end{align*}
			The transformed system reads
			\begin{align*}
				\begin{aligned}
					\dot{x}^1&=\frac{x^1}{x^2-\bar{u}^1}\,\bar{u}^2\\
					\dot{x}^2&=\bar{u}^2\\
					\dot{x}^3&=\sqrt{\frac{x^1}{x^2-\bar{u}^1}}\,\bar{u}^2
				\end{aligned}
			\end{align*}
			and according to Corollary \ref{cor:1}, we can render this system static feedback linearizable by $\#R-n=3$ prolongations of the new input $\bar{u}^1$, \ie the prolonged system
			\begin{align*}
				\begin{aligned}
					\dot{x}^1&=\frac{x^1}{x^2-\bar{u}^1}\,\bar{u}^2\,,&&&\dot{\bar{u}}^1&=\bar{u}^1_1\\
					\dot{x}^2&=\bar{u}^2\,,&&&\dot{\bar{u}}^1_1&=\bar{u}^1_2\\
					\dot{x}^3&=\sqrt{\frac{x^1}{x^2-\bar{u}^1}}\,\bar{u}^2\,,&&&\dot{\bar{u}}^1_2&=\bar{u}^1_3
				\end{aligned}
			\end{align*}
			with the new input $(\bar{u}^2,\bar{u}^1_3)$, is static feedback linearizable. 
			
			\subsubsection{Example 3.}\label{ex:3}In the following we derive a tracking control law for the planar VTOL aircraft, given by
			\begin{align*}
				\begin{aligned}
					\dot{x}&=v_x\,,&&&\dot{v}_x&=\epsilon\Cos(\theta)u^2-\Sin(\theta)u^1\\
					\dot{z}&=v_z\,,&&&\dot{v}_z&=\Cos(\theta)u^1+\epsilon\Sin(\theta)u^2-1\\
					\dot{\theta}&=\omega\,,&&&\dot{\omega}&=u^2\,.
				\end{aligned}
			\end{align*}
			It possesses the $x$-flat output $\bar{y}=(x-\epsilon\Sin(\theta),z+\epsilon\Cos(\theta))$. Applying our method to derive a tracking control law for this flat output is of course possible, but for demonstration purposes, let us instead derive a tracking control law for the $(x,u)$-flat output
			\begin{align*}
				\begin{aligned}
					y&=(x-\epsilon\Sin(\theta)+\Cos(\theta)u^1-1-\epsilon\omega^2\Cos(\theta),z+\epsilon\Cos(\theta))\,,
				\end{aligned}
			\end{align*}
			which is obtained from the $x$-flat output $\bar{y}$ by adding the second time derivative of the second component to the first one \ie $y^1=\bar{y}^1+\bar{y}^2_2$, $y^2=\bar{y}^2$. The parameterizing map with respect to this flat output is of the form
			\begin{align*}
				\begin{aligned}
					x&=F_x(y^1_{[3]},y^2_{[5]})\\
					u&=F_u(y^1_{[4]},y^2_{[6]})\,,
				\end{aligned}
			\end{align*}
			\ie $R=(4,6)$ and one easily verifies that we have $K=(0,2)$. For the input transformation \eqref{eq:6} we can choose
			\begin{align}\label{eq:25}
				\begin{aligned}
					\bar{u}^1&=x-\epsilon\Sin(\theta)+\Cos(\theta)u^1-1-\epsilon\omega^2\Cos(\theta)\\
					\bar{u}^2&=u^2\,.
				\end{aligned}
			\end{align}
			According to Theorem \ref{thm:2}, $\tilde{x}_B=y_{[\kappa-1]}$ with $\kappa=(0,6)$, \ie $\tilde{x}_B=y^2_{[5]}$, is a generalized Brunovsky state of the system, and this generalized Brunovsky state satisfies that $\partial_{\tilde{x}_B}F_x$ is regular. The regularity of $\partial_{\tilde{x}_B}F_x$ guarantees that this generalized Brunovsky state can be expressed as a function of $x$ and $y_{[\kappa,R-1]}=y^1_{[3]}=v^1_{[3]}$, which follows as
			\begin{align*}
				\begin{aligned}
					\tilde{x}_B&=\tilde{\phi}(x,v^1_{[3]})=(\tilde{x}_B^0,\ldots,\tilde{x}_B^5)\,,
				\end{aligned}
			\end{align*}
			with
			\begin{align*}
				\begin{aligned}
					\tilde{x}_B^0&=z+\epsilon\Cos(\theta)\\
					\tilde{x}_B^1&=v_z-\epsilon\omega\Sin(\theta)\\
					\tilde{x}_B^2&=v^1-x+\epsilon\Sin(\theta)\\
					\tilde{x}_B^3&=v^1_1-v_x+\epsilon\omega\Cos(\theta)\\
					\tilde{x}_B^4&=v^1_2-\epsilon\Cos(\theta)+\frac{1}{\Cos(\theta)}(\epsilon+(v^1-x+1)\Sin(\theta))\\
					\tilde{x}_B^5&=v^1_3+\epsilon\omega\Sin(\theta)+(v^1_1-v_x)\Tan(\theta)+\\				&\hspace{8ex}\frac{\omega}{c^2(\theta)}(v^1-x+1+\epsilon\Sin(\theta))
				\end{aligned}
			\end{align*}
			by applying the input transformation \eqref{eq:25}, computing the needed time derivatives $y_{[\kappa-1]}$ of the flat output\linebreak\\\\\\\\
		\begin{minipage}{\textwidth}
			\setcounter{section}{2}
			\setcounter{equation}{0}
			\begin{align}\label{eq:27}
				\begin{aligned}
					&\partial_{(x,\bar{u}^1_{[r_1-k_1-1]})}\hat{\bar{F}}_{e,sub}=&\begin{bmatrix}
						\partial_x\bar{\varphi}_{[K-1]}&0&0&0&\ldots&\ldots&\ldots&0\\
						0&1&0&0&\ldots&\ldots&\ldots&0\\
						\partial_x\bar{\varphi}^2_{k_2}&\partial_{\bar{u}^1}\bar{\varphi}^2_{k_2}&0&0&\ldots&\ldots&\ldots&0\\
						0&0&1&0&\ldots&\ldots&\ldots&0\\
						\partial_x\bar{\varphi}^2_{k_2+1}&\partial_{\bar{u}^1}\bar{\varphi}^2_{k_2+1}&\partial_{\bar{u}^1_1}\bar{\varphi}^2_{k_2+1}&0&\ldots&\ldots&\ldots&0\\
						\vdots&&&&&&&\vdots\\
						0&0&0&0&\ldots&0&1&0\\
						\partial_x\bar{\varphi}^2_{r_2-2}&\partial_{\bar{u}^1}\bar{\varphi}^2_{r_2-2}&\partial_{\bar{u}^1_1}\bar{\varphi}^2_{r_2-2}&\partial_{\bar{u}^1_2}\bar{\varphi}^2_{r_2-2}&\ldots&\ldots&\partial_{\bar{u}^1_{r_1-k_1-2}}\bar{\varphi}^2_{r_2-2}&0\\
						0&0&0&0&\ldots&\ldots&0&1\\
						\partial_x\bar{\varphi}^2_{r_2-1}&\partial_{\bar{u}^1}\bar{\varphi}^2_{r_2-1}&\partial_{\bar{u}^1_1}\bar{\varphi}^2_{r_2-1}&\partial_{\bar{u}^1_2}\bar{\varphi}^2_{r_2-1}&\ldots&\ldots&\ldots&\partial_{\bar{u}^1_{r_1-k_1-1}}\bar{\varphi}^2_{r_2-1}
					\end{bmatrix}
				\end{aligned}
			\end{align}
			\setcounter{section}{1}
			\setcounter{equation}{2}
		\end{minipage}
		\newpage
		in these coordinates and replacing the occurring time derivatives $\bar{u}^1_{[r_2-k_2-1]}$ by $v^1_{[r_2-k_2-1]}$. Together with $y_{[\kappa,R]}=(y^1_{[4]},y^2_6)=(v^1_{[4]},v^2)$, we can replace $y_{[R]}$ in $u=F_u(y_{[R]})$ to obtain a linearizing feedback of the form
		\begin{align}\label{eq:26}
			\begin{aligned}
				u&=F_u(\tilde{\phi}(x,v^1_{[3]}),v^1_{[4]},v^2)\,.
			\end{aligned}
		\end{align}
		What is left to do is to solve the equation system \eqref{eq:23}, in this example given by
			\begin{align*}
				\begin{aligned}
					v^1&=y^{1,d}\,,&&&v^1_1&=y^{1,d}_1\\
					v^1_2&=y^{1,d}_2\,,&&&v^1_3&=y^{1,d}_3\\
					v^1_4&=y^{1,d}_4\,,&&&v^2&=y^{2,d}_6-\sum_{\beta=0}^{5}a^{2,\beta}\left(\tilde{x}_B^\beta-y^{2,d}_\beta\right)\,.
				\end{aligned}
			\end{align*}
			Inserting its solution into \eqref{eq:26} gives a control law of the form
			\begin{align*}
				\begin{aligned}
					u&=\alpha(x,y^{1,d}_{[4]},y^{2,d}_{[6]})\,.
				\end{aligned}
			\end{align*}
			This control law results in the tracking error dynamics
			\begin{align*}
				\begin{aligned}
					e^1&=0\,,&&&e^2_6+\sum_{\beta=0}^{5}a^{2,\beta}e^2_\beta&=0\,.
				\end{aligned}
			\end{align*}
		\section{Proof of Theorem \ref{thm:2}}\label{app:b}
			\setcounter{equation}{1}
			In Section \ref{se:main}, we only provided a part of the proof of Theorem \ref{thm:2}. Here we complete the proof by showing that $\partial_x\bar{\varphi}_{[\kappa-1]}$ with $\kappa=(k_1,r_2)$ is regular and thus $\tilde{x}_B=y_{[\kappa-1]}$ is indeed a valid generalized Brunovsky state of the $(x,u)$-flat two-input system \eqref{eq:1}. The Jacobian matrix of \eqref{eq:11} is given by \eqref{eq:27} below. Since it is the Jacobian matrix of a diffeomorphism, its rows (as well as its columns) are linearly independent. From the rows of \eqref{eq:27} we can construct the matrix
			\begin{align}\label{eq:28}
				\begin{aligned}
					M&=\begin{bmatrix}
						\partial_x\bar{\varphi}_{[\kappa-1]}&0
					\end{bmatrix}\,,
				\end{aligned}
			\end{align}
			\ie the rows of $M$ are linear combinations of the rows of \eqref{eq:27}. Each row of $M$ is constructed by taking one line of \eqref{eq:27} corresponding to one of the $n$ components of $\bar{\varphi}_{[\kappa-1]}$ and combining it with the rows corresponding to the components $\bar{\varphi}_{[\kappa,R-1]}$. The linear independence of the rows of \eqref{eq:27} implies the linear independence of the such constructed $n$ rows of $M$, \ie $\Rank{M}=n$. The block of zeros besides the Jacobian matrix $\partial_x\bar{\varphi}_{[\kappa-1]}$ in \eqref{eq:28} does not contribute to the rank of $M$. Thus, the $n$ rows of the $n\times n$ Jacobian matrix $\partial_x\bar{\varphi}_{[\kappa-1]}$ are linearly independent, which means $\partial_x\bar{\varphi}_{[\kappa-1]}$ is regular.
\end{document}